\documentclass[a4paper,9pt]{article}
\usepackage{amsfonts}
\usepackage{amssymb,amsthm,color,tikz}
\usepackage{amsmath}
\usepackage{comment}

\usepackage[normalem]{ulem}

\newcommand{\R}{\mathbb{R}}
\newcommand{\N}{\mathbb{N}}

\newcommand{\Z}{\mathbb{Z}}

\newcommand{\E}{\mathbb{E}}
\newcommand{\Pa}{\mathbb{P}}

\newtheorem{theorem}{Theorem}
\newtheorem{lemma}[theorem]{Lemma}

\newtheorem{definition}{Definition}

\newcommand{\Om}{\Omega}

\newcommand{\bp}{\begin{proof}}
\newcommand{\ep}{\end{proof}}

\begin{document}
\title{Localisation for the torsion function and the strong Hardy inequality}

\author{{M. van den Berg} \\
School of Mathematics, University of Bristol\\
Fry Building, Woodland Road\\
Bristol BS8 1UG, United Kingdom\\
\texttt{mamvdb@bristol.ac.uk}\\
\\
{T. Kappeler}\\
Institut f\"ur Mathematik, Universit\"at Z\"urich\\
 Winterthurerstrasse 190,
CH-8057 Z\"urich, Switzerland\\
 \texttt{thomas.kappeler@math.uzh.ch}}
\date{23 February 2020}\maketitle
\vskip 3truecm \indent
\begin{abstract}\noindent
Two-sided bounds for the efficiency of the torsion function are obtained in terms of the square of the distance to the boundary function under the hypothesis that the Dirichlet Laplacian satisfies a strong Hardy inequality.
Localisation properties of the torsion function are obtained under that hypothesis. An example is analysed in detail.

\end{abstract}
\begin{center} \textbf{Keywords}: Torsion function, Dirichlet boundary condition, localisation, efficiency, Hardy inequality.
\end{center}


\section{Introduction and main results \label{sec1}}
Let $\Omega$ be an open set in $\R^m$, $m\ge 1$, with finite Lebesgue measure $|\Om|$, $0 < |\Omega| < \infty$, and with boundary $\partial\Omega$.
The torsion function for $\Omega$ is the unique solution of
\begin{equation*}
-\Delta v=1,\, \qquad v\in H_0^{ 1 }(\Omega),
\end{equation*}
and is denoted by $v_{\Omega}$.
 The function $v_{\Omega}$ is non-negative and satisfies,
\begin{equation}\label{e2}
\lambda_1(\Omega)^{-1}\le \|v_{\Omega}\|_{\infty}\le
(4+3m\log 2)\lambda_1(\Omega)^{-1},
\end{equation}
where
\begin{equation}\label{e3}
\lambda_1(\Omega)=\inf_{\varphi\in
H_0^1(\Omega)\setminus\{0\}}\frac{\Vert \nabla\varphi\Vert_2^2}{
\Vert \varphi\Vert_2^2}
\end{equation}
is the first eigenvalue of the Dirichlet Laplacian
and where $\Vert\cdot\Vert_p$ denotes the standard $L^p$ norm, $1\le p\le \infty$.
The $m$-dependent constant in the right-hand side of \eqref{e2} has
subsequently been improved (\cite{GS}, \cite{HV}). We denote the sharp constant by $\mathfrak{c}_m$,
\begin{equation}\label{e4}
\mathfrak{c}_m=\sup\{\lambda_1(\Omega)\|v_{\Omega}\|_{\infty}:\Omega\,\textup{ open in }\,\R^m,\,0<|\Omega|<\infty\}.
\end{equation}

The torsion function and its $L^1$ norm, the torsional rigidity,
play key roles in different parts of
analysis. For example the torsional rigidity of a cross section of a
beam appears in the computation of the angular change when a beam of
a given length and a given modulus of rigidity is exposed to a
twisting moment (\cite{Bandle},\cite{PSZ}). It also arises in the
definition of gamma convergence \cite{BB} and in the study of
minimal submanifolds \cite{MP}. Moreover, $v_{\Om}(x)$ equals the expected lifetime of Brownian motion in $\Om$ starting at $x\in\Om$.
This immediately implies the non-negativity of $v_{\Om}$, and that for open sets $\Om_1,\Om_2$,
\begin{equation}\label{e3a}
\Om_1\subset \Om_2 \Rightarrow v_{\Omega_1}(x)\le v_{\Omega_2}(x),\, \forall x\in \Omega_1.
\end{equation}

This in turn implies that the torsional rigidity  $\|v_{\Om}\|_1$ is monotone increasing in $\Om$.
The torsion function has been studied extensively,
and numerous works have been written on this subject. We just mention
the paper \cite{LB}, and the references therein.

The torsion function is defined  for an open subset $\Omega \subset \R^m$ of infinite measure,
provided the bottom of the spectrum of the Dirichlet Laplacian, denoted by $\lambda_1(\Omega)$, is bounded away from $0$.
Indeed, by considering an increasing sequence of sets $\Omega\cap B(0;k),\,k\in \N$,  where $B(p;R)=\{x\in\R^m:|p-x|<R\}$
denotes the open ball with centre $p$ and radius $R$, one obtains a weak solution of
\begin{equation}\label{e3b}
-\Delta v_{\Om}=1,\qquad v|_{\partial\Om} = 0.
\end{equation}
This solution is non-negative, satisfies both \eqref{e3a} and \eqref{e2}, and satisfies the probabilistic interpretation mentioned above.
See Theorem 1 in \cite {vdBC}, Theorem 5.3 in \cite{vdB}, and \cite{BB} for further details.

We remark that the torsion function has been defined in greater generality. Instead of the Dirichlet Laplacian one considers
a positive, self-adjoint Schr\"odinger operator $L:= - \Delta + V$, where $V$ is bounded, measurable, and non-negative
(or even more generally, a positive, self-adjoint, elliptic operator of second order)
acting in $L^2(\Om)$ with Dirichlet or Neumann boundary conditions or on closed manifolds. In the latter situation one requires potentials $V$ such that the bottom of the spectrum of the Laplacian is bounded away from $0$.
We denote by $v_L$ the torsion function associated to $L$. It was discovered in recent papers \cite{DDJM1}, \cite{DDJM2}, and references therein,
that under appropriate conditions, the reciprocal $v_L^{-1}$ of the torsion function $v_L$ can be used
for approximating eigenvalues and eigenfunctions of the Schr\"odinger operator $L$.
The phenomenon of localisation of eigenfunctions of Schr\"odinger operators is a prominent and very active research area
and has important applications in the applied sciences. The literature is extensive. See for example the review paper of \cite{grebenkov}. Explicit estimates on the localisation of the eigenfunctions in terms of bounds
on their exponential decay away from the subdomains, where they concentrate,  were obtained in \cite{DDJM2}.

In this paper we consider the {\em $\kappa$-localisation},  the {\em localisation}, and the {\em efficiency}
of the torsion function of the Dirichlet Laplacian on an open set $\Omega$ in $\R^m$ with $0<|\Omega|<\infty$ -- see Definition \ref{def0} and Definition \ref{def1} below.
These notions have been first introduced for Dirichlet eigenfunctions of Schr\"odinger operators and go back to \cite{grebenkov} and respectively, \cite{ps}, \cite {sperb}.
The notion of efficiency, defined as the mean to max ratio, can be viewed as a rough measure of localisation.

\smallskip

In  \cite{vdBBK}, we studied the efficiency and localisation for the torsion function of a Schr\"odinger operator $-\Delta+V$ acting in $L^2(\Om)$ with Dirichlet boundary condition.
Among other results it was shown, in Theorem 4, that under appropriate conditions localisation for the torsion function implies localisation for the first Dirichlet eigenfunction.
The converse does not hold. Consider for example a sequence of ellipsoids $(\Om_n)$ with semi-axes of length $1$ and $n$ respectively. We see by Theorem 1(iii) below that $(v_{\Om_n})$ does not localise. On the other hand it was shown in Example 10 in \cite{vdBdBDPG} that the corresponding sequence of first Dirichlet eigenfunctions localises.
Unlike the torsion function the first Dirichlet eigenfunction is not monotone on set inclusion, and this in general complicates its analysis.

In this paper, we continue this study. Our results are obtained under the hypothesis that the Dirichlet Laplacian satisfies the strong Hardy inequality, defined as follows.
\begin{definition}\label{def2}
The Dirichlet Laplacian $-\Delta$ acting in $L^2(\Omega)$ satisfies the strong Hardy
inequality, with constant $c_{\Omega}\in(0,\infty)$, if
\begin{equation}\label{e9}
\|\nabla w\|_2^2 \ge \frac{1}{c_{\Omega}} \int_{\Omega}\frac{w^2}{d_{\Omega}^2},
\quad \forall\, w \in C_c^\infty(\Omega),
\end{equation}
where $d_{\Omega}$ is the distance to the boundary function,
\begin{equation*}
d_{\Omega}(x)=\inf\{|x-y|:y\in \R^m\setminus\Omega\},\qquad x\in\Omega.
\end{equation*}
\end{definition}
Both the validity and applications
of inequalities like \eqref{e9} to spectral theory and partial differential equations have been
investigated in depth. See for example \cite{Anc},  \cite{EBD1}, \cite{EBD2}, \cite{EBD3} and \cite{EBD4}.
In particular it was shown in \cite[p. 208]{Anc}, that for any proper simply connected open subset $\Omega$ in $\R^2$,
$c_{\Omega}=16$.

\smallskip
Our first result in Theorem \ref{the1}, stated below, states that the efficiency of the torsion function of the Dirichlet Laplacian on an open set $\Omega$ of $\mathbb R^m$ with $0 < |\Omega| < \infty$
can be bounded from above and below by the efficiency of the square of the distance to the boundary function
whereas item (ii) and (iii) of Theorem \ref{the1} give sufficient conditions for the  localisation of the torsion function to hold and respectively, not to hold.

In Theorem \ref{the2}(i)--(iii) we analyse an example illustrating the use of Theorem \ref{the1}.
In Theorem \ref{the2}(iv) we give an example of $\kappa$-localisation with $0<\kappa<1$.

\bigskip

We now define the notion of efficiency and localisation in precise terms.
The efficiency of $v_{\Omega}$, introduced in \cite{DPGGLB}, is defined as follows.

\begin{definition}\label{def0}
Let $\Omega$ be an open set in $\R^m$
with $0 < |\Omega | < \infty$. The efficiency, or mean to max ratio, of $v_{\Omega}$ is
\begin{equation*}
\Phi(\Omega)=\frac{\|v_{\Omega}\|_1}{|\Omega|
\|v_{\Omega}\|_\infty}.
\end{equation*}
Let $(\Omega_n)$ be a sequence of open sets in $\R^m$ with $0<|\Omega_n|<\infty,\,n\in\N$. We say that $(v_{\Omega_n})$ has vanishing efficiency if
$\lim_{n\rightarrow\infty}\Phi(\Omega_n)=0$.
\end{definition}

The following notion of localisation for the torsion function has been motivated by the one for eigenfunctions in \cite{grebenkov}, and has been used previously in \cite{vdBBK}.

\begin{definition}\label{def1}
For any sequence of open sets $(\Omega_n)$ in $\R^m$ with $0 < |\Omega_n|<\infty,\,n\in \N$, let
\begin{equation}\label{e7}
\mathfrak{A}((\Omega_n))=\bigg\{(A_n): (\forall n\in \N)(A_n\subset\Omega_n, A_n\, \textup{measurable}), \lim_{n\rightarrow\infty}\frac{|A_n|}{|\Omega_n|}=0\bigg\},
\end{equation}
and
\begin{equation}\label{e8}
\kappa=\sup\bigg\{\limsup_{n\rightarrow\infty}\frac{\int_{A_{n}}v_{\Omega_n}}{\|v_{\Omega_n}\|_1}:(A_n)\in\mathfrak{A}((\Omega_n))\bigg\}.
\end{equation}
We say that if \textup{(i)} $0<\kappa<1$ then $(v_{\Omega_n})$ $\kappa$-{\it localises},
\textup{(ii)} $\kappa=1$ then $(v_{\Omega_n})$ {\it localises}, (iii) $\kappa=0$ then $(v_{\Omega_n})$ {\it  does not localise}.
\end{definition}

Before stating our main results we review some basic facts. We note that
\begin{equation}\label{e10a}
B(x;d_{\Omega}(x))\subset \Omega,
\end{equation}
and $|B(x;d_{\Omega}(x))|\le |\Omega|$. Hence $d_{\Omega}$ is bounded from above and
\begin{equation*}
\|d_{\Omega}\|_{\infty}\le \big(\omega_m^{-1}|\Omega|\big)^{1/m},
\end{equation*}
where $\omega_m=|B(0;1)|$. Since the torsion function is pointwise increasing with respect to the domain,\eqref{e3a}, it follows by \eqref{e10a} that,
\begin{equation}\label{e10c}
v_{\Omega}(x)\ge v_{B(x;d_{\Omega}(x))}(x)=\frac{d_{\Omega}(x)^2}{2m},\quad\forall x\in \Omega,
\end{equation}
where we have used that
\begin{equation*}
v_{B(p;R)}(x)=\frac{1}{2m}(R^2-|x-p|^2),\quad \forall x\in B(p;R).
\end{equation*}
By \eqref{e10c},
\begin{equation}\label{e24}
\frac{1}{2m}\|d_{\Omega}^2\|_1\le \|v_{\Omega}\|_1 \ .
\end{equation}
Under the additional assumption that $\Omega$ satisfies \eqref{e9},
 it was shown in \cite[Theorem 2]{vdBC} that
\begin{equation}\label{e11}
\|v_{\Omega}\|_1\le c_{\Omega}\|d_{\Omega}^2\|_1.
\end{equation}
Furthermore $d_\Omega$ is uniformly continuous. Hence for any $\eta \in \R^+,$
the subset $\{ d_\Omega \ge \eta  \}$ is relatively closed in $\Omega$, and hence measurable.

\medskip
Throughout we denote for $\nu\ge0$, by  $j_{\nu}$ the first positive zero of the Bessel function $J_{\nu}$. For the ball $B(p;1)\subset \R^m, m\ge 2$, we have
$\lambda_1(B(p;1))=j^2_{(m-2)/2}$.
Our main results are as follows.

\begin{theorem}\label{the1}
\begin{enumerate}
\item[\textup{(i)}]If $\Omega$ is an open subset of $\R^m$ with $0 < |\Omega|<\infty$, and which satisfies \eqref{e9} with the strong Hardy constant $c_{\Omega}$, then
\begin{equation}\label{e12}
(2m\mathfrak{c}_mc_{\Omega})^{-1}\frac{\Vert d^2_{\Omega}\Vert_1}{|\Omega|\Vert d^2_{\Omega}\Vert_{\infty}}\le \Phi(\Omega)\le c_{\Omega}j_{(m-2)/2}^{2}\frac{\Vert d^2_{\Omega}\Vert_1}{|\Omega|\Vert d^2_{\Omega}\Vert_{\infty}},
\end{equation}
where
$\mathfrak{c}_m$ is the sharp constant in \eqref{e4}.
\item[\textup{(ii)}]Let $(\Omega_n)$ be a sequence of open sets in $\R^m$ with $0<|\Omega_n|<\infty,\,n\in\N$, and which satisfies \eqref{e9} with strong Hardy constants $c_{\Omega_n}$.
Suppose
\begin{equation}\label{e17}
c=\sup\{c_{\Omega_n}:n\in\N\}<\infty.
\end{equation}
If  $(\eta_n)$ is a sequence of strictly positive real numbers such that
\begin{equation}\label{e18}
\lim_{n\rightarrow\infty}\frac{|\{d_{\Omega_n}\ge \eta_n\}|}{|\Omega_n|}=0,
\end{equation}
and
\begin{equation}\label{e19}
\lim_{n\rightarrow\infty}\frac{\eta_n^2|\Omega_n|}{\int_{\{d_{\Omega_n}\ge \eta_n\}}d_{\Omega_n}^2}=0,
\end{equation}
then $(v_{\Omega_n})$ localises along the sequence $A_n=\{x\in\Omega_n:d_{\Omega_n}\ge \eta_n\}$.
\item[\textup{(iii)}]Let $(\Omega_n)$ be a sequence of open sets in $\R^m$ with $0<|\Omega_n|<\infty,\,n\in\N$,
which satisfies \eqref{e9} with strong Hardy constants $c_{\Omega_n}.$
Suppose \eqref{e17} holds.
If any sequence $(A_n)$ of measurable sets, $A_n\subset\Omega_n,\,n\in \N,$ with
\begin{equation*}
\lim_{n\rightarrow\infty}\frac{|A_n|}{|\Omega_n|}=0,
\end{equation*}
satisfies
\begin{equation*}
\lim_{n\rightarrow\infty}\frac{\int_{A_n}d_{\Omega_n}^2}{\int_{\Omega_n}d_{\Omega_n}^2}=0,
\end{equation*}
then $(v_{\Omega_n})$ does not localise.
\end{enumerate}
\end{theorem}

Theorem \ref{the1}(i) can be interpreted as follows. Given an open subset $\Omega$ of $\R^m$ with $0 < |\Omega|<\infty,$ we
define the efficiency of $d_\Omega^2$ as
$$
D(\Omega) = \frac{\| d_\Omega^2 \|_1}{|\Omega| \| d_\Omega^2 \|_\infty} \ .
$$
Theorem \ref{the1}(i) asserts that under condition \eqref{e9}, the efficiencies of $v_\Omega$ and $d_\Omega^2$ are comparable.

Let $(\Omega_n)$ be a sequence of open sets in $\R^m$ with $0 < |\Omega_n|<\infty$.
We say that $(d_{\Omega_n}^2)$ localises if
$$
\sup\bigg\{\limsup_{n\rightarrow\infty}\frac{\int_{A_{n}}d^2_{\Omega_n}}{\|d^2_{\Omega_n}\|_1}:(A_n)\in\mathfrak{A}((\Omega_n))\bigg\} = 1,
$$
where $\mathfrak{A}((\Omega_n)) $
is given by \eqref{e7}.
Theorem \ref{the1}(iii) asserts that for any sequence $(\Omega_n)$ satisfying the conditions of Theorem \ref{the1}(iii)
the following holds:  if $(d_{\Omega_n}^2)$ does not localise, then neither does $(v_{\Omega_n})$.

\medskip

In Theorem \ref{the2} below we analyse the torsion function for a sequence of open simply connected sets $(\Omega_{\varepsilon_n,n})_{n\ge 4}\subset \R^2$ with $0<\varepsilon_n<1$. The construction is as follows.
Let $Q$ be the open unit square in $\R^2$ with vertices $(0,0),(1,0),(0,1),(1,1)$. Let $n\in \N$, and for any given $n\in \N$ let $L_1,...,L_{n-1}$ be the closed line segments of
lengths $1-\varepsilon_n$ with endpoints $(\frac1n,0),(\frac2n,0),...,(\frac{n-1}{n},0)$ pointing in the direction $(0,1)$.  Let
$$\Omega_{\varepsilon_n,n}=Q\setminus\cup _{j=1}^{n-1}L_j.$$

\medskip

\begin{figure}[!ht]
\centering
\begin{tikzpicture}
\draw[black, thick] (0,0) rectangle (8,8);
\draw[black, thick](0.5,0) -- (0.5,7.50);
\draw[black, thick](1,0) -- (1,7.50);
\draw[black, thick](1.5,0) -- (1.5,7.50);
\draw[black, thick](2,0) -- (2,7.50);
\draw[black, thick](2.5,0) -- (2.5,7.50);
\draw[black, thick](3,0) -- (3,7.50);
\draw[black, thick](3.5,0) -- (3.5,7.50);
\draw[black, thick](4,0) -- (4,7.50);
\draw[black, thick](4.5,0) -- (4.5,7.50);
\draw[black, thick](5,0) -- (5,7.50);
\draw[black, thick](5.5,0) -- (5.5,7.50);
\draw[black, thick](6,0) -- (6,7.50);
\draw[black, thick](6.5,0) -- (6.5,7.50);
\draw[black, thick](7,0) -- (7,7.50);
\draw[black, thick](7.5,0) -- (7.5,7.50);

\end{tikzpicture}
\caption{$\Omega_{\frac{1}{16},16}$}
\label{fig1}
\end{figure}
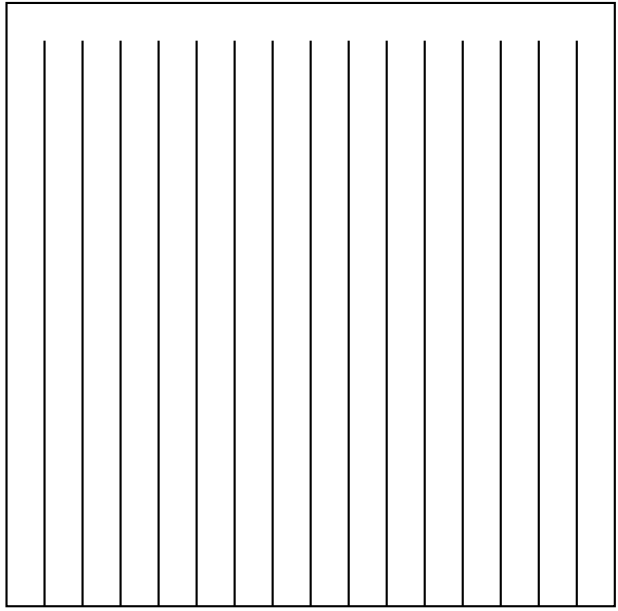

\medskip

\begin{theorem}\label{the2}
Let $0<\alpha<1,\, c>0$,
\begin{equation*}
\varepsilon_n=cn^{-\alpha},
\end{equation*}
and
\begin{equation*}
N_{\alpha,c}=\min\{n\in\N:n^{\alpha}\ge 2c,\, cn^{1-\alpha}\ge 2\}.
\end{equation*}
\textup{(}Note that $N_{\alpha,c}\ge 4$.\textup{)}
\begin{itemize}
\item[\textup{(i)}] If $n\ge N_{\alpha,c}$, then
\begin{equation}\label{e40}
\frac{1}{3072\mathfrak{c}_2}\big(cn^{-\alpha}+c^{-2}n^{2\alpha-2}\big)\le \Phi(\Omega_{\varepsilon_n,n})\le\frac{64j_0^2}{3}\big(cn^{-\alpha}+c^{-2}n^{2\alpha-2}\big),
\end{equation}
where $\mathfrak{c}_2$ is given by \eqref{e4}.
\item[\textup{(ii)}]If $0<\alpha<\frac23$, then $(v_{\Omega_{\varepsilon_n,n}})$ localises along the sequence
\begin{equation}\label{e40a}
A_n=\{x\in\Omega_{\varepsilon_n,n}:d_{\Omega_{\varepsilon_n,n}}(x)\ge \frac{1}{2n}\}.
\end{equation}
\item[\textup{(iii)}]If $\frac23<\alpha<1$, then $(v_{\Omega_{\varepsilon_n,n}})$ is not localising.
\item[\textup{(iv)}]If $\alpha=\frac23$, then $(v_{\Omega_{\varepsilon_n,n}})$ is $\kappa_c$-localising along \eqref{e40a} with
\begin{equation}\label{e41}
\kappa_c=\frac{c^3}{1+c^3}.
\end{equation}
\end{itemize}
\end{theorem}

\medskip

The paper is organised as follows. The proofs of Theorem \ref{the1} and Theorem \ref{the2}(i)--(iii) are given in Section \ref{sec2}.
The proof of Theorem \ref{the2}(iv) involves tools from Brownian motion and is deferred to Section \ref{sec3}.

\section{Proofs of Theorem \ref{the1} and Theorem \ref{the2}(i)--(iii) \label{sec2}}

{\it Proof of Theorem \textup{\ref{the1}(i).}}
Since $|\Omega|<\infty$, the inradius $\|d_{\Omega}\|_{\infty}$ is finite, and by \eqref{e3} and \eqref{e9},
\begin{align}\label{e13}
\lambda_1(\Omega)&\ge c_{\Omega}^{-1}\inf_{\varphi\in
H_0^1(\Omega),\,\|\varphi\|_2^2=1}\int_{\Omega}\frac{\varphi^2}{d_{\Omega}^2}\nonumber \\ &
\ge c_{\Omega}^{-1}\|d_{\Omega}\|_{\infty}^{-2}\inf_{\varphi\in
H_0^1(\Omega),\,\|\varphi\|_2^2=1}\Vert \varphi\Vert_2^2\nonumber \\ &
=c_{\Omega}^{-1}\|d_{\Omega}\|_{\infty}^{-2}.
\end{align}
By \eqref{e4} and \eqref{e13},
\begin{equation}\label{e14}
\|v_{\Omega}\|_{\infty}\le \mathfrak{c}_m\lambda_1(\Omega)^{-1}\le \mathfrak{c}_mc_{\Omega}\|d_{\Omega}\|_{\infty}^{2}.
\end{equation}
The lower bound in \eqref{e12} follows from \eqref{e14} and the lower bound in \eqref{e24}.

Since $\Omega$ has inradius $\|d_{\Omega}\|_{\infty}$, $\Omega$ contains an open ball with radius $\|d_{\Omega}\|_{\infty}$ and centre $p_{\Om}$. By the monotonicity of the Dirichlet eigenvalues
\begin{equation}\label{e15}
\lambda_1(\Omega)\le \lambda_1(B(p_{\Om};\|d_{\Omega}\|_{\infty})=j_{(m-2)/2}^2\|d_{\Omega}\|_{\infty}^{-2}.
\end{equation}
By the first inequality in \eqref{e2}, and \eqref{e15},
\begin{equation}\label{e16}
\|v_{\Omega}\|_{\infty}^{-1}\le j_{(m-2)/2}^{2}\|d_{\Omega}\|_{\infty}^{-2}.
\end{equation}
The upper bound in \eqref{e12} follows from \eqref{e16} and the upper bound in \eqref{e11}.

\medskip
\noindent
{\it Proof of Theorem \textup{\ref{the1}(ii).}} By \eqref{e9}, and suppressing $n$ dependence,
\begin{align}\label{e20}
\int_{\Omega}v_{\Omega}&=-\int_{\Omega}v_{\Omega}\Delta v_{\Omega}=\int_{\Omega}|\nabla v_{\Omega}|^2\nonumber \\ &
\ge \frac{1}{c_{\Omega}}\int_{\Omega}\frac{v_{\Omega}^2}{d_{\Omega}^2}\ge \frac{1}{c_{\Omega}}\int_{\{d_{\Omega}<\eta\}}\frac{v_{\Omega}^2}{d_{\Omega}^2}\nonumber \\ &
\ge \frac{1}{c_{\Omega}\eta^2}\int_{\{d_{\Omega}<\eta\}}v_{\Omega}^2.
\end{align}
On the other hand by Cauchy-Schwarz,
\begin{equation}\label{e21}
\bigg(\int_{\{d_{\Omega}<\eta\}}v_{\Omega}\bigg)^2\le |\{d_{\Omega}<\eta\}|\int_{\{d_{\Omega}<\eta\}}v_{\Omega}^2.
\end{equation}
Combining \eqref{e20} and \eqref{e21} yields,
\begin{equation}\label{e22}
\bigg(\int_{\{d_{\Omega}<\eta\}}v_{\Omega}\bigg)^2\le c_{\Omega}\eta^2|\Omega|\bigg(\int_{\{d_{\Omega}<\eta\}}v_{\Omega}+\int_{\{d_{\Omega}\ge\eta\}}v_{\Omega}\bigg).
\end{equation}
Solving the quadratic inequality \eqref{e22} gives,
\begin{align}\label{e23}
\int_{\{d_{\Omega}<\eta\}}v_{\Omega}&\le \frac{c_{\Omega}\eta^2|\Omega|}{2}+\bigg(\frac{c_{\Omega}^2\eta^4|\Omega|^2}{4}+c_{\Omega}\eta^2|\Omega|\int_{\{d_{\Omega}\ge\eta\}}v_{\Omega}\bigg)^{1/2}\nonumber \\ &
\le c_{\Omega}\eta^2|\Omega|+\bigg(c_{\Omega}\eta^2|\Omega|\int_{\{d_{\Omega}\ge\eta\}}v_{\Omega}\bigg)^{1/2}.
\end{align}
By \eqref{e10c},
\begin{equation}\label{e25}
\int_{\{d_{\Omega}\ge\eta\}}v_{\Omega}\ge \frac{1}{2m}\int_{\{d_{\Omega}\ge \eta\}}d_{\Omega}^2.
\end{equation}
By \eqref{e23} and \eqref{e25}
\begin{equation}\label{e26}
\frac{\int_{\{d_{\Omega}<\eta\}}v_{\Omega}}{\int_{\{d_{\Omega}\ge\eta\}}v_{\Omega}}\le \frac{2mc_{\Omega}\eta^2|\Omega|}{\int_{\{d_{\Omega}\ge\eta\}}d_{\Omega}^2}+\bigg(\frac{2mc_{\Omega}\eta^2|\Omega|}{\int_{\{d_{\Omega}\ge\eta\}}d_{\Omega}^2}\bigg)^{1/2}.
\end{equation}

To complete the proof we suppose that $(\Omega_n)$ and $(\eta_n)$ are sequences satisfying \eqref{e18} and \eqref{e19} respectively, and that $c_{\Omega_n}$ satisfies \eqref{e17}. If
\begin{equation*}
\lim_{n\rightarrow\infty}\frac{\eta_n^2|\Omega_n|}{\int_{\{d_{\Omega_n}\ge\eta_n\}}d_{\Omega_n}^2}=0,
\end{equation*}
then by \eqref{e26}
\begin{equation*}
\lim_{n\rightarrow\infty}\frac{\int_{\{d_{\Omega_n}<\eta_n\}}v_{\Omega_n}}{\int_{\{d_{\Omega_n}\ge\eta_n\}}v_{\Omega_n}}=0,
\end{equation*}
and
\begin{equation*}
\lim_{n\rightarrow\infty}\frac{\int_{\{d_{\Omega_n}<\eta_n\}}v_{\Omega_n}}{\int_{\Omega_n}v_{\Omega_n}}=0.
\end{equation*}
Hence
\begin{equation}\label{e30}
\lim_{n\rightarrow\infty}\frac{\int_{\{d_{\Omega_n}\ge\eta_n\}}v_{\Omega_n}}{\int_{\Omega_n}v_{\Omega_n}}=1.
\end{equation}
Let $A_n=\{d_{\Omega_n}\ge \eta_n\}$. Then $A_n,\,n\in\N,$ is relatively closed in $\Omega$, and hence measurable. By hypothesis \eqref{e18}, $$\lim_{n\rightarrow\infty}\frac{|A_n|}{|\Omega_n|}=0,$$ and by \eqref{e30} we have $\limsup_{n\rightarrow\infty}\frac{\int_{A_{n}}v_{\Omega_n}}{\|v_{\Omega_n}\|_1}=1$. Hence \eqref{e8} holds with $\kappa=1$.

\medskip
\noindent
{\it Proof of Theorem \textup{\ref{the1}(iii).}}
By Cauchy-Schwarz, and suppressing $n$-dependence,
\begin{equation}\label{e33}
\bigg(\int_{A}v_{\Omega}\bigg)^2\le \bigg(\int_{A}d_{\Omega}^2\bigg)\bigg(\int_A\frac{v_{\Omega}^2}{d_{\Omega}^2}\bigg).
\end{equation}
By \eqref{e9}, the first two equalities in \eqref{e20}, and by \eqref{e33} we obtain
\begin{align*}
\int_{\Omega}v_{\Omega}&\ge \frac{1}{c_{\Omega}}\int_{A}\frac{v_{\Omega}^2}{d_{\Omega}^2}\nonumber \\ &
\ge \frac{1}{c_{\Omega}}\frac{\big(\int_{A}v_{\Omega}\big)^2}{\int_{A}d_{\Omega}^2}.
\end{align*}
Hence
\begin{equation*}
\int_{A}v_{\Omega}\le c_{\Omega}^{1/2}\bigg(\int_{\Omega}v_{\Omega}\bigg)^{1/2}\bigg(\int_{A}d_{\Omega}^2\bigg)^{1/2},
\end{equation*}
and together with \eqref{e24}
\begin{align*}
\frac{\int_A v_{\Omega}}{\int_{\Omega}v_{\Omega}}\le c_{\Omega}^{1/2}\frac{\bigg(\int_A d_{\Omega}^2\bigg)^{1/2}}{\bigg(\int_{\Omega}v_{\Omega}\bigg)^{1/2}}
\le (2mc_{\Omega})^{1/2}\bigg(\frac{\int_A d_{\Omega}^2}{\int_{\Omega}d_{\Omega}^2}\bigg)^{1/2}.
\end{align*}
This implies the assertion of Theorem \ref{the1}(iii).
\hfill$\square$

\medskip
\noindent
{\it Proof of Theorem \textup{\ref{the2}(i).}}
We use Theorem \ref{the1}(i). First note that for a rectangle $R_{a,b}$ with side lengths $a$ and $b$, where $a\ge b$ one has
\begin{equation}\label{e42}
\|d_{R_{a,b}}\|_2^2=\frac{ab^3}{12}-\frac{b^4}{24}\ge\frac{ab^3}{24}.
\end{equation}
To obtain a lower bound for $\Vert d^2_{\Omega_{\varepsilon_n,n}}\Vert_1$ we add a closed line segment of length $1$ to the boundary of $\Omega_{\varepsilon_n,n}$, and obtain that the resulting set is the union of $n+1$ disjoint rectangles, denoted by $\tilde{\Omega}_{\varepsilon_n,n}$.
See Figure 2.
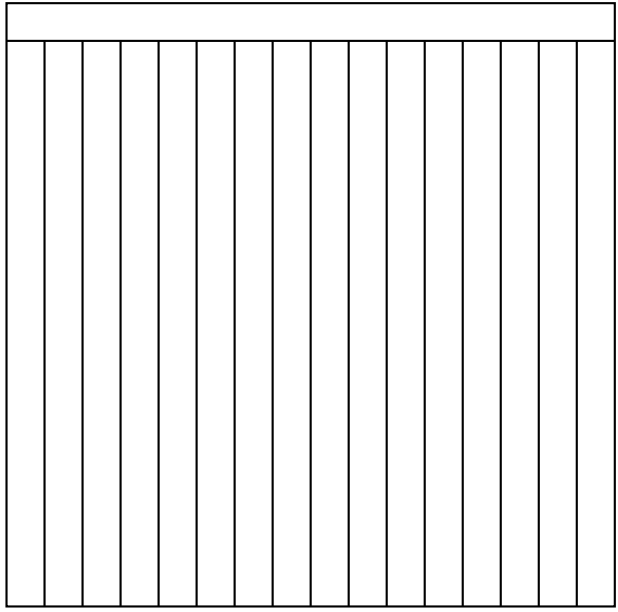
\begin{figure}[!ht]
\centering
\begin{tikzpicture}
\draw[black, thick] (0,0) rectangle (8,8);
\draw[black, thick](0.5,0) -- (0.5,7.5);
\draw[black, thick](1,0) -- (1,7.5);
\draw[black, thick](1.5,0) -- (1.5,7.5);
\draw[black, thick](2,0) -- (2,7.5);
\draw[black, thick](2.5,0) -- (2.5,7.5);
\draw[black, thick](3,0) -- (3,7.5);
\draw[black, thick](3.5,0) -- (3.5,7.5);
\draw[black, thick](4,0) -- (4,7.5);
\draw[black, thick](4.5,0) -- (4.5,7.5);
\draw[black, thick](5,0) -- (5,7.5);
\draw[black, thick](5.5,0) -- (5.5,7.5);
\draw[black, thick](6,0) -- (6,7.5);
\draw[black, thick](6.5,0) -- (6.5,7.5);
\draw[black, thick](7,0) -- (7,7.5);
\draw[black, thick](7.5,0) -- (7.5,7.5);
\draw[black, thick](0,7.5) -- (8,7.5);
\end{tikzpicture}
\caption{$\tilde{\Omega}_{\frac{1}{16},16}$}
\label{fig2}
\end{figure}
Hence, by \eqref{e42},
\begin{equation}\label{e43}
\|d^2_{\Omega_{\varepsilon_n,n}}\|_1\ge \|d^2_{\tilde{\Omega}_{\varepsilon_n,n}}\|_1\ge \frac{1-\varepsilon_n}{24n^2}+\frac{\varepsilon_n^3}{24}\ge \frac{1}{48}(\varepsilon_n^3+n^{-2}),
\end{equation}
where we have used that $\varepsilon_n\le \frac12$ for $n\ge N_{\alpha,c}$.

To obtain an upper bound for $\|d^2_{\Omega_{\varepsilon_n,n}}\|_1$, we have the following estimate.
\begin{equation}\label{e51}
d_{\Omega_{\varepsilon_n,n}}(x)\le \begin{cases}\frac{1}{2n}, \qquad \forall(x_1,x_2)\in \{\Omega_{\varepsilon_n,n}:x_2\le 1-\varepsilon_n\},\\
1-x_2,\qquad \forall(x_1,x_2)\in \{\Omega_{\varepsilon_n,n}:x_2\ge1-\varepsilon_n\}.\end{cases}
\end{equation}
Since
\begin{equation}\label{e51a}
|\{\Omega_{\varepsilon_n,n}:x_2\le 1-\varepsilon_n\}|=1-\varepsilon_n\le 1,
\end{equation}
we have by \eqref{e51} and \eqref{e51a},
\begin{equation}\label{e51b}
\int_{\{\Omega_{\varepsilon_n,n}:x_2\le 1-\varepsilon_n\}}d^2_{\Omega_{\varepsilon_n,n}}\le \frac{1}{4n^2}
\end{equation}
and
\begin{equation}\label{e51c}
\int_{\{\Omega_{\varepsilon_n,n}:x_2\ge 1-\varepsilon_n\}}d^2_{\Omega_{\varepsilon_n,n}}\le \int_{1-\varepsilon_n}^1dx_2(1-x_2)^2 =\frac{\varepsilon_n^3}{3}.
\end{equation}
Hence we have
\begin{equation}\label{e44}
\|d^2_{\Omega_{\varepsilon_n,n}}\|_1\le \frac{1}{3}(\varepsilon_n^3+n^{-2}).
\end{equation}
Finally we observe that
\begin{equation}\label{e45}
\frac12\varepsilon_n\le \Vert d_{\Omega_{\varepsilon_n,n}}\Vert_{\infty}\le \varepsilon_n,\qquad n\ge N_{\alpha,c},
\end{equation}
and \eqref{e40} follows from \eqref{e43}, \eqref{e44}, \eqref{e45}, Theorem \ref{the1}(i), $|\Omega_{\varepsilon_n,n}|=1$, and that for a proper simply connected subset $\Omega$ in $\R^2$, $c_{\Omega}=16$, \cite{Anc}. This proves the assertion of Theorem \ref{the2}(i).

\medskip
\noindent{\it Proof of Theorem \textup{\ref{the2}(ii).}} We use Theorem \ref{the1}(ii), and choose $\eta_n=\frac{1}{2n}$. Then
\begin{equation*}
|\{d_{\Omega_{\varepsilon_n,n}}\ge \frac{1}{2n}\}|\le \varepsilon_n,
\end{equation*}
and \eqref{e18} holds.
Since for $n\ge N_{\alpha,c}$, $\{d_{\Omega_{\varepsilon_n,n}}\ge \frac{1}{2n}\}$ contains a rectangle with side lengths $1-n^{-1}$ and $cn^{-\alpha}-\frac{1}{n}$ respectively, we have by \eqref{e42}
\begin{align*}
\int_{\{d_{\Omega_{\varepsilon_n,n}}\ge \frac{1}{2n}\}}d_{\Omega_{\varepsilon_n,n}}^2&\ge \frac{1}{24}(1-n^{-1})\big(cn^{-\alpha}-n^{-1}\big)^3\nonumber \\&
\ge \frac{c^3}{192}(1-n^{-1})n^{-3\alpha}\nonumber \\ &
\ge \frac{c^3}{256}n^{-3\alpha},
\end{align*}
where we have used that $cn^{1-\alpha}\ge 2$ for $n\ge N_{\alpha,c}\ge 4$.
Hence
\begin{equation*}
\frac{\eta_n^2|\Omega_{\varepsilon_n,n}|}{\int_{\{d_{\Omega_{\varepsilon_n,n}}\ge \eta_n\}}d_{\Omega_{\varepsilon_n,n}}^2}\le \frac{64n^{3\alpha-2}}{c^3},\qquad n\ge N_{\alpha,c},
\end{equation*}
and \eqref{e19} holds for $0<\alpha<\frac23$.
This proves the assertion of Theorem \ref{the2}(ii) by Theorem \ref{the1}(ii).

\medskip
\noindent{\it Proof of Theorem \textup{\ref{the2}(iii).}} We use Theorem \ref{the1}(iii), and let for $n\in \N$, $A_n$ be an arbitrary measurable subset of $\Omega_{\varepsilon_n,n}$ such that $\lim_{n\rightarrow\infty}|A_n|=0$.
By \eqref{e43},
\begin{equation}\label{e49}
\|d^2_{\Omega_{\varepsilon_n,n}}\|_1\ge \frac{1}{48n^2}.
\end{equation}
By \eqref{e51},
\begin{equation}\label{e52}
\int_{A_n\cap\{x_2\le1-\varepsilon_n\}}d_{\Omega_{\varepsilon_n,n}}^2\le \frac{|A_n|}{4n^2}.
\end{equation}
By \eqref{e51c}, \eqref{e49}, \eqref{e52},
\begin{equation*}
\frac{\int_{A_n}d_{\Omega_{\varepsilon_n,n}}^2}{\int_{\Omega_{\varepsilon_n,n}}d_{\Omega_{\varepsilon_n,n}}^2}\le  c^3 n^{2-3\alpha}+12|A_n|.
\end{equation*}
This implies
\begin{equation*}
\lim_{n\rightarrow\infty}\frac{\int_{A_n}d_{\Omega_{\varepsilon_n,n}}^2}{\int_{\Omega_{\varepsilon_n,n}}d_{\Omega_{\varepsilon_n,n}}^2}=0,
\end{equation*}
since $\frac23<\alpha<1$, and the hypothesis on $A_n$. This proves the assertion of Theorem \ref{the2}(iii) by Theorem \ref{the1}(iii).
\hfill$\square$

\section{Proof of Theorem \ref{the2}(iv)} \label{sec3}

Let $(B(s),s\ge 0; \Pa_x,x\in \R^d)$ be Brownian motion associated to the Laplacian in $\R^d$.
Here $(B(s),s\ge 0)$ takes values in $\R^d$, $\Pa_x$ is Wiener measure with $\Pa_x(B(0)=x)=1$, and for every Borel set $A\in \R^d$,
\begin{equation}\label{e54a}
\Pa_x(B(s)\in A)=(4\pi s)^{-d/2}\int_Ady\,e^{-|x-y|^2/(4s)},\quad s>0.
\end{equation}
If $\Omega\subset \R^d$ is open and $x\in \Omega$ then we denote the first exit time (or life time) of
$\Omega$ by
\begin{equation*}
T_{\Omega}=\inf\{s\ge 0: B(s)\in \R^2\setminus \Omega\}.
\end{equation*}
It is convenient to denote the first hitting time of a closed set $C$ by
\begin{equation*}\label{e56}
\tau_C=\inf\{s\ge 0: B(s)\in C\}.
\end{equation*}
It is well-known that the unique weak solution $u_{\Omega}$ of
\begin{equation*}
\Delta u=\frac{\partial u}{\partial t},\qquad \textup{in}\, \Omega\times(0,\infty),
\end {equation*}
with $u|_{\partial \Omega \times (0, \infty)} =0 $,
and with initial datum
\begin{equation*}
u(x;0)=1,\qquad x\in\Omega,
\end{equation*}
is given by
\begin{equation}\label{e59}
u_{\Omega}(x;t)=\Pa_x(T_{\Omega}>t).
\end{equation}
Since
\begin{equation}\label{e60}
v_{\Omega}(x)=\int_0^{\infty}dt\, u_{\Omega}(x;t)
\end{equation}
we obtain by \eqref{e59} and \eqref{e60} that
\begin{equation*}
v_{\Omega}(x)=\E_x(T_{\Omega}),
\end{equation*}
and
\begin{equation*}
\|v_{\Omega}\|_1=\int_{\Omega}dx\,\int_0^{\infty}dt\, u_{\Omega}(x;t).
\end{equation*}
For the facts above we refer to \cite{Dur}.
Now let $d=2$, and let $B=(B_1,B_2),$ where $B_1$ and $B_2$ are independent one-dimensional Brownian motions with probability measures $\Pa_{x_1}^{(1)}$, and $\Pa_{x_2}^{(1)}$ respectively, and $\Pa_x=\Pa_{x_1}^{(1)}\otimes\Pa_{x_2}^{(2)}$.
 We have that $\inf\{s\ge 0:B_2(s)=a\}=T_{(-\infty,a)}$. Since
\begin{equation*}
\Pa^{(2)}_{x_2}[T_{(-\infty,a)}>t]=\frac{2}{\pi^{1/2}}\int_0^{(a-x_2)/(4t)^{1/2}}d\theta e^{-\theta^2},\qquad \forall x_2<a,
\end{equation*}
we have that the density $\rho(a,\tau)$ of the random variable $\inf\{s\ge 0:B_2(s)=a\}$ with $x_2=0$, is given by
\begin{equation}\label{e64}
\rho(a,\tau)=\frac{a}{2\pi^{1/2}\tau^{3/2}}e^{-a^2/(4\tau)}{\bf1}_{\R^+}(\tau).
\end{equation}

\medskip

The following lemma will be used in the proof of  Theorem \ref{the2}(iv).

\begin{lemma}\label{lem3}
Let $\Omega$ be an open, bounded and connected set in $\R^2$ which contains an open rectangle $R_{a,b}\equiv R_{a,b}(p),\,a\ge b,$ with sides
$K_1=[(p,0)-(p+b,0)],\, K_2=[(p,0)-(p,a)],\,K_3=[(p+b,0)-(p+b,a)],\, K_4=[(p+b,a)-(p,a)]$.
Let
\begin{equation*}
\Omega_{a,b}\equiv \Omega_{a,b}(p)=\Omega\setminus\big(\cup_{j=1}^3K_j\big).
\end{equation*}
If $x\in R_{a,b}$, then
\begin{equation}\label{e65}
v_{\Omega_{a,b}}(x)\le \frac12(x_1-p)(p+b-x_1)+2^{9/2}e^{-\pi(a-x_2)/(2b)}\lambda_1(\Omega_{a,b})^{-1}.
\end{equation}
\end{lemma}

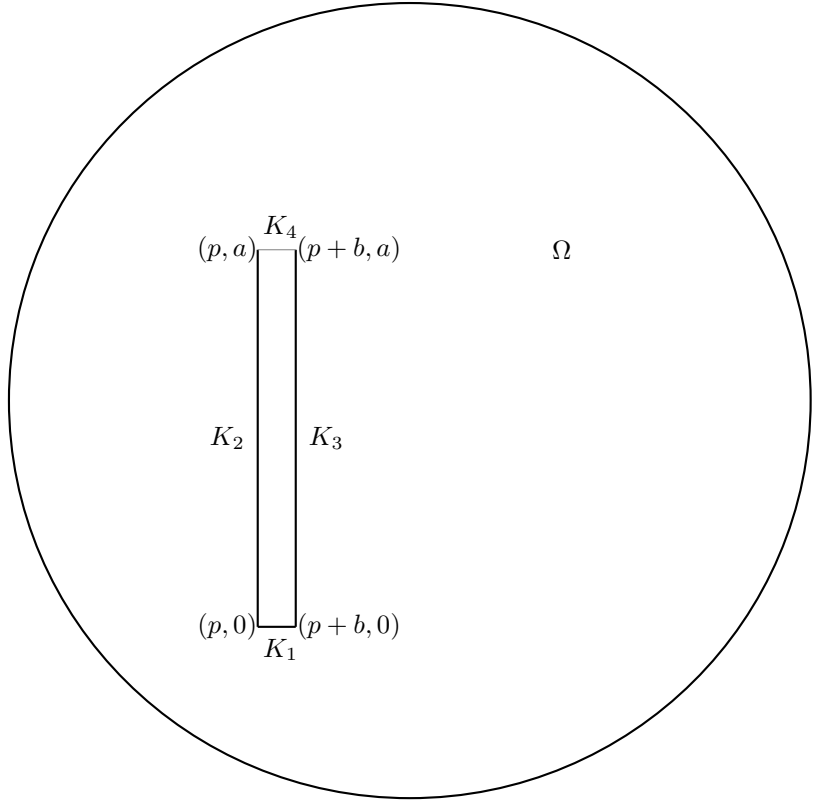
\begin{figure}[!ht]
\centering
\begin{tikzpicture}
\draw[black, thick] (6,4) circle (150pt);
\draw[black, thick](4.0,1.0) -- (4.0,6.0);
\draw[black, thick](4.5,1.0) -- (4.5,6.0);
\draw[black, thick](4.0,1.0) -- (4.5,1.0);
\draw[gray, ultra thin](4.0,6.0) -- (4.5,6.0);
\node at (3.6,3.5)    {$K_2$};
\node at (4.3,0.7)    {$K_1$};
\node at (4.9,3.5)    {$K_3$};
\node at (4.3,6.3)    {$K_4$};
\node at (8.0,6.0)    {$\Omega$};
\node at (3.6,1.0)    {$(p,0)$};
\node at (3.6,6.0)    {$(p,a)$};
\node at (5.2,1.0)    {$(p+b,0)$};
\node at (5.2,6.0)    {$(p+b,a)$};
\end{tikzpicture}
\caption{$\Omega_{a,b}=\Omega\setminus\big(\cup_{j=1}^3K_j\big)$:\quad Rectangle $R_{a,b}$
with sides $K_2, K_3$ of length $a$, and $K_1, K_4$ of length $b$.}
\label{fig3}
\end{figure}

\begin{proof}
By \eqref{e59} and \eqref{e60}, for any $x\in R_{a,b}$,
\begin{align}\label{e66}
v_{\Omega_{a,b}}(x)&=\int_0^{\infty}dt\,\Pa_x[T_{\Omega_{a,b}}>t]\nonumber \\ &
=\int_0^{\infty}dt\,\Pa_{x}[T_{R_{a,b}}>t]+\int_0^{\infty}dt\,\Pa_{x}[T_{R_{a,b}}\le t< T_{\Omega_{a,b}}]\nonumber \\ &
\le v_{R_{a,b}}(x_1,x_2)+\int_0^{\infty}dt\,\Pa_{x}[T_{R_{a,b}}\le t< T_{\Omega_{a,b}}]\nonumber \\ &
\le \frac12(x_1-p)(p+b-x_1)+\int_0^{\infty}dt\,\Pa_{x}[T_{R_{a,b}}\le t< T_{\Omega_{a,b}}],
\end{align}
where we have used that $a\mapsto v_{R_{a,b}}(x_1,x_2)$ is monotone increasing and bounded from above by $\frac12(x_1-p)(p+b-x_1)$. See Figure \ref{fig3}. The latter is the torsion function for the interval $(p,p+b)$.
To bound the second term in the right-hand side of \eqref{e66} we use the inclusion,
\begin{align*}
\{T_{R_{a,b}}\le t< T_{\Omega_{a,b}}\}\subset \{B(T_{R_{a,b}})\in K_4\}\cap\{t< T_{\Omega_{a,b}}\}.
\end{align*}
This inclusion states that the event of a Brownian path starting
for example at $x$ in $R_{a,b}$
exiting the rectangle $R_{a,b}$ before $t$ but not exiting $\Omega_{a,b}$ before $t$ has to exit the rectangle at $K_4$, while staying in $\Omega_{a,b}$ until $t$.
By Cauchy-Schwarz,
\begin{align}\label{e68}
\Pa_x[T_{R_{a,b}}\le t< T_{\Omega_{a,b}}]&\le \big(\Pa_x[B(T_{R_{a,b}})\in K_4]\big)^{1/2}(\Pa_x[t< T_{\Omega_{a,b}}])^{1/2}\nonumber \\ &
\le 2^{3/4}\big(\Pa_x[B(T_{R_{a,b}})\in K_4]\big)^{1/2}e^{-t\lambda_1(\Omega_{a,b})/8},
\end{align}
where we have used that for open sets $\Omega\subset \R^2$, or for open subsets $\Omega\subset\R^1$,
\begin{equation} \label{e68a}
\Pa_x[T_{\Om}>t]\le 2^{3/2}e^{-t\lambda_1(\Om)/4}.
\end{equation}
See pp.10, 11 in \cite{vdBBK}. Let $S_{a,b}=(p,p+b)\times(-\infty,a)$ be the half-strip of width $b$ below the line segment $K_4$.
Note that $R_{a,b} \subset S_{a,b}$ and that $|S_{a,b}| = \infty.$
By \eqref{e68}, by applying \eqref{e68a} to the open set $\Omega=(0,b)$, and using \eqref{e64},
\begin{align}\label{e69}
\Pa_x[B(T_{R_{a,b}})\in K_4]&\le\Pa_x[B(T_{S_{a,b}})\in K_4]\nonumber \\ &
=\int_0^{\infty}d\tau\,\rho(a-x_2,\tau)\Pa^{(1)}_{x_1}[T_{(0,b)}>\tau]\nonumber \\ &
\le 2^{3/2}\int_0^{\infty}d\tau\,\frac{a-x_2}{2\pi^{1/2}\tau^{3/2}}e^{-(a-x_2)^2/(4\tau)-\pi^2\tau/(4b^2)}\nonumber \\ &
=2^{3/2}e^{-\pi(a-x_2)/(2b)},
\end{align}
and where we have used formula 3.472.3 in \cite{RG} for the last equality.
By \eqref{e68} and \eqref{e69}
\begin{equation}\label{e70}
\Pa_x[T_{R_{a,b}}\le t< T_{\Omega_{a,b}}]\le 2^{3/2}e^{-\pi(a-x_2)/(4b)-t\lambda_1(\Omega_{a,b})/8}.
\end{equation}
Integrating both sides of \eqref{e70} with respect to $t$, and using \eqref{e66} completes the proof.
\end{proof}
\noindent{\it Proof of Theorem \textup{\ref{the2}(iv)}}. To prove \eqref{e41} we first obtain a lower bound on the right-hand side of \eqref{e8} by making a particular choice for $(A_n)$:
\begin{equation}\label{e71}
A_n=\{x\in\Omega_{\varepsilon_n,n}:x_2>1-cn^{-2/3}\}.
\end{equation}
We let
\begin{equation}\label{e72}
B_n=\{x\in\Omega_{\varepsilon_n,n}:1-cn^{-2/3}\ge x_2>1-cn^{-2/3}-n^{-3/4}\},
\end{equation}
and
\begin{equation}\label{e72a}
C_n=\{x\in\Omega_{\varepsilon_n,n}:1-cn^{-2/3}-n^{-3/4}\ge x_2\},
\end{equation}
so that $\Omega_{\varepsilon_n,n}=A_n\cup B_n\cup C_n$. By \eqref{e8}, and the monotonicity of the torsion function,
\begin{align}\label{e73}
\kappa&\ge \limsup_{n\rightarrow\infty}\frac{\int_{A_n}v_{\Omega_{\varepsilon_n,n}}}{\int_{A_n}v_{\Omega_{\varepsilon_n,n}}+\int_{B_n}v_{\Omega_{\varepsilon_n,n}}+\int_{C_n}v_{\Omega_{\varepsilon_n,n}}}\nonumber \\ &
\ge \limsup_{n\rightarrow\infty}\frac{\int_{A_n}v_{A_n}}{\int_{A_n}v_{A_n}+\int_{B_n}v_{\Omega_{\varepsilon_n,n}}+\int_{C_n}v_{\Omega_{\varepsilon_n,n}}}.
\end{align}
Since $A_n$ is a rectangle with side lengths $1$ and $cn^{-2/3}$ respectively, we have by \eqref{e42}
\begin{equation}\label{e74}
\frac{c^3}{12n^2}\ge \int_{A_n}v_{A_n}\ge \frac{c^3}{12n^2}-\frac{c^4}{24n^{8/3}}.
\end{equation}
By \eqref{e73} and \eqref{e74},
\begin{align}\label{e75}
\kappa&\ge \limsup_{n\rightarrow\infty}\frac{c^3-\frac12c^4n^{-2/3}}{c^3+12n^2\int_{B_n}v_{\Omega_{\varepsilon_n,n}}+12n^2\int_{C_n}v_{\Omega_{\varepsilon_n,n}}}\nonumber \\ &
=\limsup_{n\rightarrow\infty}\frac{c^3}{c^3+12n^2\int_{B_n}v_{\Omega_{\varepsilon_n,n}}+12n^2\int_{C_n}v_{\Omega_{\varepsilon_n,n}}}
\end{align}
By \eqref{e14} and \eqref{e45},
\begin{equation}\label{e76}
\|v_{\Omega_{\varepsilon_n,n}}\|_{\infty}\le 16\mathfrak{c}_2\|d_{\Omega_{\varepsilon_n,n}}\|^2_{\infty}\le 16\mathfrak{c}_2c^2n^{-4/3},\qquad n\ge N_{2/3,c}.
\end{equation}
Since $|B_n|=n^{-3/4}$ we have by \eqref{e76},
\begin{align*}
12n^2\int_{B_n}v_{\Omega_{\varepsilon_n,n}}&\le  12n^2|B_n| \|v_{\Omega_{\varepsilon_n,n}}\|_{\infty} \nonumber \\ &\le 192\mathfrak{c}_2c^2n^{-1/12}.
\end{align*}
This, together with \eqref{e75}, yields
\begin{equation}\label{e78}
\kappa\ge \limsup_{n\rightarrow\infty}\frac{c^3}{c^3+12n^2\int_{C_n}v_{\Omega_{\varepsilon_n,n}}}.
\end{equation}
In order to bound the integral in the right-hand side from above we use Lemma \ref{lem3} for each of the $n$ rectangles $R_{a,b}(p)$ in $\Omega_{\varepsilon_n,n}$, with $p=\frac{k}{n},\,0\le k\le n-1,$
with $b=\frac1n$, and $a=1-cn^{-2/3}$. Note that for any $p=\frac{k}{n}$ the set $\Omega_{a,b}(p)$ introduced in Lemma \ref{lem3} coincides with $\Omega_{\varepsilon_n,n}$, and that each point in $C_n$
satisfies $a-x_2\ge n^{-3/4}$. By \eqref{e13} and \eqref{e45},
\begin{equation*}
\lambda_1(\Omega_{a,b})=\lambda_1(\Omega_{\varepsilon_n,n})\ge \frac{1}{16\|d_{\Omega_{\varepsilon_n,n}}\|^2_{\infty}}\ge \frac{1}{16c^2n^{-4/3}}.
\end{equation*}
Hence the second term the right-hand side of \eqref{e65} is bounded from above by
$2^{17/2}c^2n^{-4/3}e^{-\pi n^{1/4}/2}$
uniformly for all points $x\in C_n$. Since $|C_n|\le 1$ we have
\begin{equation}\label{e80}
\int_{C_n}2^{9/2}e^{-\pi(a-x_2)/(2b)}\lambda_1(\Omega_{a,b})^{-1}\le  2^{17/2}c^2n^{-4/3}e^{-\pi n^{1/4}/2}.
\end{equation}
Integrating the first term in the right-hand side of \eqref{e65} over one rectangle in $C_n,$ and adjusting the coordinate frame appropriately, gives a contribution
\begin{equation*}
\int_0^{1-cn^{-2/3}-n^{-3/4}}dx_2\int_0^{1/n}dx_1\frac{x_1}{2}(n^{-1}-x_1)=\frac{1}{12n^3}\big(1-cn^{-2/3}-n^{-3/4}\big).
\end{equation*}
Summing over all $n$ rectangles in $C_n$ gives, together with \eqref{e80} and Lemma \ref{lem3},
\begin{equation}\label{e82}
12n^2\int_{C_n}v_{\Omega_{\varepsilon_n,n}}\le 1-cn^{-2/3}-n^{-3/4}+3\cdot2^{21/2}c^2n^{2/3}e^{-\pi n^{1/4}/2}.
\end{equation}
By \eqref{e78} and \eqref{e82} we conclude
\begin{equation*}
\kappa\ge \frac{c^3}{1+c^3}.
\end{equation*}

We now prove the converse inequality by obtaining an upper bound for the right-hand side of \eqref{e8}. We first observe, by the monotonicity, that the torsion for $\Omega_{\varepsilon_n,n}$ is bounded from below by the torsion for the set in Figure \ref{fig2}. The latter is the union of $n$ disjoint rectangles $R_{a,b}$ with $b=\frac1n$ and $a=1-cn^{-2/3}$, together with one disjoint rectangle with $a=1$ and $b=cn^{-2/3}$.
By \eqref{e42} we find
\begin{equation}\label{e84}
\int_{\Omega_{\varepsilon_n,n}}v_{\Omega_{\varepsilon_n,n}}\ge \frac{1-cn^{-2/3}}{12n^2}-\frac{1}{24n^3}+\frac{c^3}{12n^2}-\frac{c^4}{24n^{8/3}}.
\end{equation}
In order to avoid abuse of notation we keep \eqref{e71}, \eqref{e72} and \eqref{e72a}, and denote by $(D_n)$ an arbitrary sequence of in $\mathfrak{A}((\Omega_{\varepsilon_n,n}))$. That is $D_n\subset \Omega_{\varepsilon_n,n}$, measurable, $n\in \N,$ with
$\lim_{n\rightarrow\infty}|D_n|=0$. We have
\begin{align}\label{e85}
\int_{D_n}v_{\Omega_{\varepsilon_n,n}}&=\int_{D_n\cap(A_n\cup B_n)}v_{\Omega_{\varepsilon_n,n}}+\int_{D_n\setminus(A_n\cup B_n)}v_{\Omega_{\varepsilon_n,n}}\nonumber \\ &
\le\int_{A_n}v_{\Omega_{\varepsilon_n,n}}+\int_{B_n}v_{\Omega_{\varepsilon_n,n}}+\int_{D_n\cap C_n}v_{\Omega_{\varepsilon_n,n}}.
\end{align}
By \eqref{e76},
\begin{equation}\label{e86}
\int_{B_n}v_{\Omega_{\varepsilon_n,n}}\le 16\mathfrak{c}_2c^2n^{-25/12}.
\end{equation}
We use Lemma \ref{lem3} to bound the third term in the right-hand side of \eqref{e85} from above. We have, as before, that for each point $x\in C_n$, both $a-x_2\ge n^{-3/4}$, and
$\frac{x_1}{2}(\frac1n-x_1)\le \frac{1}{8n^2}$, where the $x_1$ coordinate is adjusted to the rectangle under consideration. This gives,
\begin{equation*}
v_{\Omega_{\varepsilon_n,n}}(x)\le \frac{1}{8n^2}+2^{17/2}c^2n^{-4/3}e^{-\pi n^{1/4}/2},\quad \forall x\in C_n.
\end{equation*}
Hence
\begin{align}\label{e88}
\int_{D_n\cap C_n}v_{\Omega_{\varepsilon_n,n}}&\le |D_n\cap C_n|\Big(\frac{1}{8n^2}+2^{17/2}c^2n^{-4/3}e^{-\pi n^{1/4}/2}\Big)\nonumber \\& \le|D_n|\Big(\frac{1}{8n^2}+2^{17/2}c^2n^{-4/3}e^{-\pi n^{1/4}/2}\Big).
\end{align}
By \eqref{e84} and \eqref{e86},
\begin{equation*}
\limsup_{n\rightarrow\infty}\frac{\int_{B_n}v_{\Omega_{\varepsilon_n,n}}}{\int_{\Omega_{\varepsilon_n,n}}v_{\Omega_{\varepsilon_n,n}}}=0.
\end{equation*}
Moreover, by \eqref{e84} and \eqref{e88},
\begin{equation*}
\limsup_{n\rightarrow\infty}\frac{\int_{D_n\cap C_n}v_{\Omega_{\varepsilon_n,n}}}{\int_{\Omega_{\varepsilon_n,n}}v_{\Omega_{\varepsilon_n,n}}}\le\limsup_{n\rightarrow\infty}\frac{4|D_n|}{3(1+c^3)}=0.
\end{equation*}

\begin{figure}[!ht]
\centering
\begin{tikzpicture}
\draw[black, thick] (-0.75,8) -- (9.75,8);
draw [black, thick] (-2.0,-1)--(-2.0,7.5);
draw [black, thick] (-1.5,-1)--(-1.5,7.5);
draw [black, thick] (-1.0,-1)--(-1.0,7.5);
draw [black, thick] (-0.5,-1)--(-0.5,7.5);
draw [black, thick] (0.0,-1)--(0.0,7.5 );
\draw[black, thick](0.5,-1) -- (0.5,7.5);
\draw[black, thick](1,-1) -- (1,7.5);
\draw[black, thick](1.5,-1) -- (1.5,7.5);
\draw[black, thick](2,-1) -- (2,7.5);
\draw[black, thick](2.5,-1) -- (2.5,7.5);
\draw[black, thick](3,-1) -- (3,7.5);
\draw[black, thick](3.5,-1) -- (3.5,7.5);
\draw[black, thick](4,-1) -- (4,7.5);
\draw[black, thick](4.5,-1) -- (4.5,7.5);
\draw[black, thick](5,-1) -- (5,7.5);
\draw[black, thick](5.5,-1) -- (5.5,7.5);
\draw[black, thick](6,-1) -- (6,7.5);
\draw[black, thick](6.5,-1) -- (6.5,7.5);
\draw[black, thick](7,-1) -- (7,7.5);
\draw[black, thick](7.5,-1) -- (7.5,7.5);
\draw[black, thick](8.0,-1) -- (8.0,7.5);
\draw[black, thick](8.5,-1) -- (8.5,7.5);
\draw[black, thick](9.0,-1) -- (9.0,7.5);
\draw[black, thick] (9.5,-1)--(9.5,7.5);
\node at (-0.0,6.0)    {$\tilde{\Omega}_n$};
\end{tikzpicture}
\caption{$\tilde{\Omega}_{n}$:\quad The vertical half lines, periodically extended, are at distance $\frac1n$, and have distance $cn^{-2/3}$ to the horizontal line.}
\label{fig4}
\end{figure}
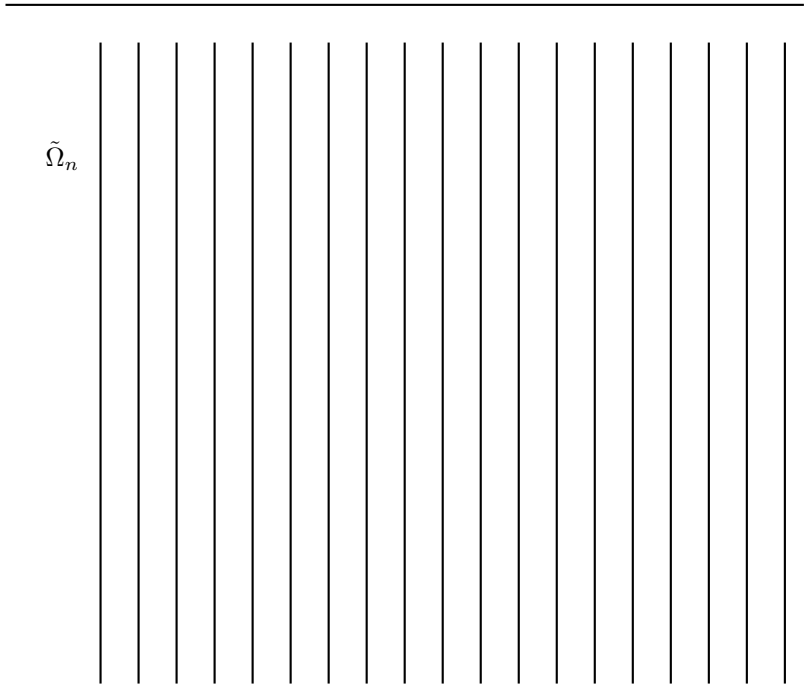

It remains to find an upper bound for $\int_{A_n}v_{\Omega_{\varepsilon_n,n}}$. Let $\tilde{\Omega}_n$ be the connected open set in $\R^2$ with boundary consisting of the horizontal line $\{x_2=1\}$ and the vertical half-lines $\{x_2\le 1-cn^{-2/3}, x_1=\frac{k}{n},\, k\in \Z\}$. See Figure \ref{fig4}.

We first show that the bottom of the spectrum of the Dirichlet Laplacian on $\tilde{\Omega}_n$ is bounded away from $0$.
By inserting Neumann boundary conditions on the line $x_2=1-cn^{-2/3}$ the bottom of the spectrum decreases and
decouples into the union of the spectra of the horizontal strip at width $cn^{-2/3}$ with Dirichlet boundary conditions on $x_2=1$,
and Neumann boundary condition on $x_2=1-cn^{-2/3}$, and vertical half-strips of width $n^{-1}$ with Neumann boundary condition
at the top end and Dirichlet conditions on the half lines. By taking the double of these half-strips
we obtain that the bottom of the spectrum here is the same as the bottom of the infinite strip of width $n^{-1}$.
Combining these inequalities gives for $n\ge N_{\alpha,c}$,
\begin{align}\label{e90a}
\lambda_1(\tilde{\Omega}_n)&\ge \min\{\frac{\pi^2n^{4/3}}{4c^2},\pi^2n^2\}\nonumber \\ &
=\frac{\pi^2n^{4/3}}{c^2}.
\end{align}
Hence the torsion function for $\tilde{\Omega}_n$ exists (see \eqref{e3b} and the preceding paragraph), and by the monotonicity property \eqref{e3a}
\begin{equation*}
v_{\Omega_{\varepsilon_n,n}}(x)\le v_{\tilde{\Omega}_n}(x),\quad \forall x\in \Omega_{\varepsilon_n,n}.
\end{equation*}
Note that $$v_{\tilde{\Omega}_n}(x_1,x_2)=v_{\tilde{\Omega}_n}(x_1+\frac1n,x_2),\quad \forall x\in \tilde{\Omega}_n.$$
Let
\begin{equation*}
E_n=\{x\in \R^2:1> x_2>1-cn^{-2/3}-n^{-3/4}\}.
\end{equation*}
We have
\begin{align}\label{e93}
\Pa_x[T_{\Omega_{\varepsilon_n,n}}>t]&\le \Pa_x[T_{\tilde{\Omega}_n}>t]\nonumber \\ &
=\Pa_x[T_{\tilde{\Omega}_n}>t,\, T_{E_n}>t]+\Pa_x[T_{E_n}\le t<T_{\tilde{\Omega}_n}]\nonumber \\ &
\le\Pa_x[T_{E_n}>t]+\Pa_x[T_{E_n}\le t<T_{\tilde{\Omega}_n}].
\end{align}
Integrating both sides of \eqref{e93} with respect to $t$ over $\R^+$ gives
\begin{align*}
v_{\Omega_{\varepsilon_n,n}}(x)&\le v_{E_n}(x)+\int_0^{\infty}dt\,\Pa_x[T_{E_n}\le t<T_{\tilde{\Omega}_n}]\nonumber \\ &
=\frac12(1-x_2)(x_2-1+d_n)+\int_0^{\infty}dt\,\Pa_x[T_{E_n}\le t<T_{\tilde{\Omega}_n}],
\end{align*}
where
\begin{equation}\label{e95}
d_n=cn^{-2/3}+n^{-3/4}.
\end{equation}
Hence
\begin{align*}
\int_{A_n}v_{\Omega_{\varepsilon_n,n}}&\le \int_{A_n}dx\,\frac12(1-x_2)(x_2-1+d_n)+\int_{A_n}dx\int_0^{\infty}dt\,\Pa_x[T_{E_n}\le t<T_{\tilde{\Omega}_n}]\nonumber \\ &
\le\frac{1}{12}d_n^3+\int_{A_n}dx\int_0^{\infty}dt\,\Pa_x[T_{E_n}\le t<T_{\tilde{\Omega}_n}].
\end{align*}
By \eqref{e84}, \eqref{e95} we have
\begin{equation*}
\limsup_{n\rightarrow\infty}\frac{\frac{1}{12}d_n^3}{\int_{\Omega_{\varepsilon_n,n}}v_{\Omega_{\varepsilon_n,n}}}\le \frac{c^3}{1+c^3}.
\end{equation*}

To complete the proof it therefore suffices to show that
\begin{equation}\label{e98}
\limsup_{n\rightarrow\infty}n^2\int_{A_n}dx\int_0^{\infty}dt\,\Pa_x[T_{E_n}\le t<T_{\tilde{\Omega}_n}]=0.
\end{equation}
To prove \eqref{e98} we let $L_n=\{(x_1,1-d_n):x_1\in \R\}$. We have by the strong Markov property,
\begin{align}\label{e99}
\Pa_x[T_{E_n}\le t<T_{\tilde{\Omega}_n}]&\le \mathbb{E}_x\Big(\int_0^t {\bf 1}_{\tau_{L_n}\in d\tau} \Pa_{B(\tau_{L_n})}[T_{\tilde{\Omega}_n}>t-\tau]\Big)\nonumber \\ &
\le \mathbb{E}_x\Big(\int_0^t {\bf 1}_{\tau_{L_n}\in d\tau} \sup_{z\in L_n}\Pa_{z}[T_{\tilde{\Omega}_n}>t-\tau]\Big)\nonumber \\ &
=\int_0^t d\tau\frac{\partial}{\partial \tau}\big(\Pa_x[\tau_{L_n}<\tau]\big) \sup_{z\in K_n}\Pa_{z}[T_{\tilde{\Omega}_n}>t-\tau]\Big)\nonumber \\ &
=\int_0^t d\tau\frac{\partial}{\partial \tau}\big(\Pa_x[\tau_{L_n}<\tau]\big)\Pa_{((2n)^{-1},1-d_n)}[T_{\tilde{\Omega}_n}>t-\tau],
\end{align}
where ${\bf1}_{\tau_{L_n}\in d\tau}$ is the indicator function on the set of Brownian paths which hit $L_n$ in the infinitesimal interval $d\tau$.
We have used in the final equality in \eqref{e99} that $x_1\mapsto \Pa_{(x_1,1-d_n)} [T_{\tilde{\Omega}_n}>t-\tau],x_1\in\R$ is periodic in $x_1$ with period $n^{-1}$, and which has equal maxima at
$\{((2n)^{-1}+kn^{-1},1-d_n):k\in \Z\}$.
Hence the supremum in the third line in \eqref{e99} is a maximum for $z=((2n)^{-1},1-d_n)$.
Integrating the convolution in the right-hand side of \eqref{e99} with respect to $t$ over $\R^+$ gives
\begin{align*}
\int_0^{\infty} dt\,\Pa_x[T_{E_n}\le\, &t<T_{\tilde{\Omega}_n}]\nonumber \\ &\le \Big(\int _0^\infty d\tau\frac{\partial}{\partial \tau}\big(\Pa_x[\tau_{L_n}<\tau]\big)\Big)\int_0^\infty d\tau\,\Pa_{((2n)^{-1},1-d_n)}[T_{\tilde{\Omega}_n}>\tau]\nonumber \\ & =v_{\tilde{\Omega}_n}((2n)^{-1},1-d_n).
\end{align*}
One verifies that Lemma \ref{lem3} also holds if $\Omega$ contains the open rectangle $(p,p+b)\times(-\infty,a)$. Indeed the first term in the right-hand side of \eqref{e65} is the torsion function for the infinite strip of width $b$. The second term takes into account that the strip is cut-off at $a$. With this we arrive at
\begin{equation}\label{e101}
\int_0^t dt\,\Pa_x[T_{E_n}\le t<T_{\tilde{\Omega}_n}]\le \frac{1}{8n^2}+2^{9/2}e^{-\pi n^{1/4}/2}\lambda_1(\tilde{\Omega}_n)^{-1}.
\end{equation}
Since $|A_n|=cn^{-2/3}$, we find by \eqref{e101}, and \eqref{e90a},
\begin{equation*}
n^2\int_{A_n}dx\int_0^t dt\,\Pa_x[T_{E_n}\le t<T_{\tilde{\Omega}_n}]\le \frac{c}{8n^{2/3}}+2^{9/2}\pi^{-2}c^3e^{-\pi n^{1/4}/2}.
\end{equation*}
This implies \eqref{e98}.
\hfill$\square$

\bigskip

\bigskip

{\bf Acknowledgments.} Both authors acknowledge support by the Leverhulme Trust through Emeritus Fellowship EM-2018-011-9,
and the Swiss National Science Foundation respectively.

\end{document}